\documentclass{amsart}

\usepackage{amsmath}
\usepackage{amsfonts}
\usepackage{amssymb}
\usepackage{graphicx}
\usepackage{hyperref}
\usepackage{mathrsfs}
\usepackage{pb-diagram}
\usepackage{epstopdf}
\usepackage{amsmath,amsfonts,amsthm,enumerate,amscd,latexsym}
\usepackage{curves}
\usepackage{bbm}
\usepackage[mathscr]{eucal}
\usepackage{caption}
\usepackage{subcaption}
\usepackage{tikz}
\newtheorem{thm}{Theorem}[section]

\newtheorem{prop}[thm]{Proposition}
\newtheorem{defn}[thm]{Definition}
\newtheorem{example}[thm]{Example}
\newtheorem{remark}[thm]{Remark}
\newtheorem{conj}[thm]{Conjecture}

\numberwithin{equation}{section}

\def\bZ{\mathbb{Z}}
\def\bQ{\mathbb{Q}}
\def\bR{\mathbb{R}}
\def\bC{\mathbb{C}}
\def\bP{\mathbb{P}}

\begin{document}

\title[SYZ for toric varieties]{SYZ mirror symmetry for toric varieties}
\author[K. Chan]{Kwokwai Chan}
\address{Department of Mathematics\\ The Chinese University of Hong Kong\\ Shatin\\ Hong Kong}
\email{kwchan@math.cuhk.edu.hk}

\date{\today}

\begin{abstract}
We survey recent developments in the study of SYZ mirror symmetry for compact toric and toric Calabi-Yau varieties, with a special emphasis on works of the author and his collaborators.
\end{abstract}

\maketitle

\tableofcontents

\section{The SYZ proposal}

In 1996, drawing on the new idea of D-branes from string theory, Strominger, Yau and Zaslow \cite{SYZ96} made a ground-breaking proposal to explain mirror symmetry geometrically as a Fourier--type transform, known as {\em T-duality}.

\begin{conj}[The SYZ conjecture \cite{SYZ96}]
If $X$ and $\check{X}$ are Calabi-Yau manifolds mirror to each other, then
\begin{itemize}
\item[(i)]
both $X$ and $\check{X}$ admit special Lagrangian torus fibrations with sections $\rho:X\to B$ and $\check{\rho}:\check{X}\to B$ over the same base which are fiberwise dual to each other in the sense that regular fibers $\rho^{-1}(b) \subset X$ and $\check{\rho}^{-1}(b) \subset \check{X}$ over the same point $b \in B$ are dual tori, and

\item[(ii)]
there exist Fourier--type transforms which exchange symplectic-geometric (resp. complex-geometric) data on $X$ with complex-geometric (resp. symplectic-geometric) data on $\check{X}$.
\end{itemize}
\end{conj}

This conjecture has since become one of the two major mathematical approaches in the study of mirror symmetry (the other being Kontsevich's homological mirror symmetry conjecture \cite{Kontsevich-ICM94}); it provides a concrete and beautiful picture depicting the geometry underlying mirror symmetry.

The SYZ conjecture, in particular, suggests a geometric construction of mirrors, namely, by taking the fiberwise dual of a special Lagrangian torus fibration on a given Calabi-Yau manifold. It turns out that such a construction has to be modified by {\em quantum (or instanton) corrections} except in the semi-flat case (which will be reviewed below).

By {\em SYZ mirror symmetry} we mean the investigation of mirror symmetry along the lines of the SYZ proposal. The aim of this note is to review recent developments in SYZ mirror symmetry for toric varieties, including SYZ constructions of their mirrors and various open mirror theorems which lead to an enumerating meaning of (inverse of) mirror maps.

\subsection{Semi-flat SYZ: a toy example}

In the so-called {\em semi-flat} case where special Lagrangian torus fibrations do not admit any singular fibers, the SYZ picture is particularly appealing. To describe this, we need the following

\begin{defn}
An {\em affine structure} on an $n$-dimensional smooth manifold is an atlas of charts such that the transition maps lie in the affine linear group $\text{Aff}(\bR^n) = \bR^n \rtimes GL_n(\bR)$. The affine structure is called {\em tropical} (resp. {\em integral}) if the transition maps lie in $\bR^n \rtimes GL_n(\bZ)$ (resp. $\text{Aff}(\bZ^n) = \bZ^n \rtimes GL_n(\bZ)$).
\end{defn}

Given a special Lagrangian torus fibration $\rho: X\to B$ {\em without} singular fibers,\footnote{In this note, a Lagrangian submanifold $L \subset X$ is called {\em special} if $\text{Im}\ \Omega|_L = 0$, where $\Omega$ is the holomorphic volume form on $X$.} classic results of McLean \cite{McLean98} give us two naturally defined tropical affine structures on the base $B$: the {\em symplectic} and {\em complex} affine structures. In these terms, mirror symmetry can be explained neatly as an interchange of these structures.

To describe these structures, let $L_b := \rho^{-1}(b)$ be a fiber of $\rho$. A tangent vector $v \in T_b B$ can be lifted to a vector field (also denoted by $v$ by abuse of notations) normal to the fiber $L_b$, and this determines the differential forms
\begin{align*}
\alpha & := -\iota_v \omega \in \Omega^1(L_b;\bR), \\
\beta & := \iota_v \text{Im}\ \Omega \in \Omega^{n-1}(L_b;\bR),
\end{align*}
where $\omega$ and $\Omega$ denote the K\"ahler form and holomorphic volume form on the Calabi-Yau manifold $X$ respectively.

McLean \cite{McLean98} proved that the corresponding deformation of $L_b$ gives special Lagrangian submanifolds if and only if both $\alpha$ and $\beta$ are closed. By identifying the tangent space $T_b B$ with the cohomology group $H^1(L_b;\bR)$ by $v \mapsto -\iota_v \omega$, we obtain the symplectic affine structure on $B$; similarly, identifying $T_b B$ with the cohomology group $H^{n-1}(L_b;\bR)$ by $v \mapsto \iota_v \text{Im}\ \Omega$ produces the complex affine structure on $B$. We also have the {\em McLean metric} on $B$ defined by
$$g(v_1,v_2) := -\int_{L_b} \iota_{v_1}\omega \wedge \iota_{v_2}\text{Im}\ \Omega.$$

In a very interesting paper \cite{Hitchin97}, Hitchin explains how these tropical affine structures are related in an elegant way by the {\em Legendre transform}. Let $x_1,\ldots,x_n$ be local affine coordinates on $B$ with respect to the symplectic affine structure, then the McLean metric can be written as the Hessian of a convex function $\phi$ on $B$, i.e.
$$g\left(\frac{\partial}{\partial x_j}, \frac{\partial}{\partial x_k}\right) = \frac{\partial^2\phi}{\partial x_j\partial x_k}.$$
Hitchin observed that setting $\check{x}_j := \partial\phi/\partial x_j$ ($j=1,\ldots,n$) gives precisely the local affine coordinates on $B$ with respect to the complex affine structure. Furthermore, let
$$\check{\phi} := \sum_{j=1}^n \check{x}_j x_j - \phi(x_1,\ldots,x_n)$$
be the Legendre transform of $\phi$. Then we have $x_j = \partial\check{\phi}/\partial\check{x}_j$ and
$$g\left(\frac{\partial}{\partial \check{x}_j}, \frac{\partial}{\partial \check{x}_k}\right) = \frac{\partial^2\check{\phi}}{\partial \check{x}_j\partial \check{x}_k}.$$

Now suppose that the fibration $\rho:X \to B$ admits a Lagrangian section. Then a theorem of Duistermaat \cite{Duistermaat80} gives global {\em action-angle coordinates} and hence a symplectomorphism
$$X \cong T^*B/\Lambda^\vee,$$
where $\Lambda^\vee \subset T^*B$ denotes the lattice locally generated by $dx_1,\ldots,dx_n$, and the quotient $T^*B/\Lambda^\vee$ is equipped with the canonical symplectic structure
$$\omega_{\text{can}} := \sum_{j=1}^n dx_j \wedge dy^j,$$
where $y^1,\ldots,y^n$ are fiber coordinates on $T^*B$.

The SYZ proposal suggests that the mirror of $X$ is the total space of the fiberwise dual of $\rho: X \to B$, which is nothing but the quotient
$$\check{X} := TB/\Lambda,$$
where $\Lambda \subset TB$ is the dual lattice locally generated by $\partial/\partial x_1,\ldots,\partial/\partial x_n$. Geometrically, $\check{X}$ should be viewed as the moduli space of pairs $(L,\nabla)$ where $L$ is a fiber of $\rho$ and $\nabla$ is a flat $U(1)$-connection over $L$; this is in fact the key idea lying at the heart of the SYZ proposal.

The quotient $\check{X} = TB/\Lambda$ is naturally a complex manifold with holomorphic coordinates given by $z_j := \exp(x_j + \mathbf{i}y_j)$, where $y_1,\ldots,y_n$ are fiber coordinates on $TB$ dual to $y^1,\ldots,y^n$. So the SYZ approach constructs the mirror of $X$ as a complex manifold equipped with a nowhere vanishing holomorphic volume form
$$\check{\Omega} := d\log z_1 \wedge \cdots \wedge d\log z_n$$
and a dual torus fibration $\check{\rho}: \check{X} \to B$.

In this case, one can write down an explicit fiberwise Fourier-type transform, which we call the {\em semi-flat SYZ transform} $\mathcal{F}^\text{semi-flat}$, that carries $e^{\mathbf{i}\omega}$ to $\check{\Omega}$. Indeed, by straightforward computations, we have
\begin{align*}
\check{\Omega} & = \int_{\rho^{-1}(b)}e^{\mathbf{i}\omega} \wedge e^{\mathbf{i}\sum_j dy_j \wedge dy^j},\\
e^{\mathbf{i}\omega} & = \int_{\check{\rho}^{-1}(b)}\check{\Omega} \wedge e^{-\mathbf{i}\sum_j dy_j \wedge dy^j}.
\end{align*}
These can naturally be interpreted as Fourier-Mukai--type transforms because the differential form $e^{\mathbf{i}\sum_j dy_j \wedge dy^j}$ is the Chern character for a universal connection on the fiberwise Poincar\'e bundle $\mathcal{P}$ on $X \times_B \check{X}$:
\begin{prop}
We have
\begin{align*}
\check{\Omega} & = \mathcal{F}^\text{semi-flat}\left(e^{\mathbf{i}\omega}\right),\\
e^{\mathbf{i}\omega} & = \left(\mathcal{F}^\text{semi-flat}\right)^{-1}\left(\check{\Omega}\right),
\end{align*}
where the semi-flat SYZ transforms are defined as
\begin{align*}
\mathcal{F}^\text{semi-flat}\left(-\right) & := \pi_{2\ast}\left(\pi_{1}^{\ast}\left(-\right)\wedge \text{ch}(\mathcal{P})\right),\\
\left(\mathcal{F}^\text{semi-flat}\right)^{-1}\left(-\right) & := \pi_{1\ast}\left(\pi_{2}^{\ast}\left(-\right)\wedge \text{ch}(\mathcal{P})^{-1}\right);
\end{align*}
here $\pi_{j}$ denotes projection to the $j^{\text{th}}$-factor.
\end{prop}

This explains one-half of mirror symmetry in the semi-flat case, namely, the correspondence between the A-model on $X$ and B-model on $\check{X}$.
If we now switch to the complex affine structure on $B$, then the symplectic structure
$$\check{\omega} := \mathbf{i}\partial\bar{\partial}\phi = \sum_{j,k} \phi_{jk} dx_j\wedge dy_k,$$
where $\phi_{jk} = \frac{\partial^2\phi}{\partial x_j\partial x_k}$, on $\check{X}$ is compatible with its complex structure, and thus $\check{X}$ becomes a K\"ahler manifold.

On the other side, we obtain a complex structure on $X$ with holomorphic coordinates $w_i,\ldots,w_n$ defined by
$d\log w_j = \sum_{k=1}^n \phi_{jk}dx_k + \mathbf{i}dy^j$,
so $X$ also becomes a K\"ahler manifold, and is equipped with the nowhere vanishing holomorphic $n$-form
$$\Omega := d\log w_1 \wedge \cdots \wedge d\log w_n.$$
The maps $\rho: X \to B$ and $\check{\rho}: \check{X} \to B$ now give fiberwise dual special Lagrangian torus fibrations.

\begin{remark}
When the function $\phi$ above satisfies the {\em real Monge-Amp\`ere equation}
$$\text{det}\left(\frac{\partial^2\phi}{\partial x_j\partial x_k}\right) = \text{constant},$$
we obtain $T^n$-invariant {\em Ricci-flat} metrics on both $X$ and its mirror $\check{X}$. In this case, the McLean metric on $B$ is called a {\em Monge-Amp\`ere metric} and $B$ will be called a {\em Monge-Amp\`ere manifold}. This explains why SYZ mirror symmetry is intimately related to the study of real Monge-Amp\`ere equations and affine K\"ahler geometry. We refer the interested readers to the papers \cite{Kontsevich-Soibelman01, Loftin05, LYZ05} and references therein for the detail story.
\end{remark}

The SYZ conjecture therefore paints a beautiful picture for mirror symmetry in the semi-flat case; see \cite{LYZ00, Leung05} for other details on semi-flat SYZ mirror symmetry (cf. also \cite[Section 2]{Chan-Leung10b} and \cite[Chapter 6]{D-branes-MS_book}). Unfortunately, this nice picture remains valid only at the large complex structure/volume limits where all quantum corrections are suppressed. Away from these limits, special Lagrangian torus fibrations will have singular fibers and the SYZ mirror construction must be modified by instanton corrections which, as suggested again by SYZ \cite{SYZ96}, should come from holomorphic disks bounded by the Lagrangian torus fibers.

\section{A quick review on toric varieties}

Before discussing mirror symmetry for toric varieties, let us set up some notations and review certain basic facts in toric geometry which we need later. We will follow standard references such as the books \cite{Fulton_toric_book, CLS_toric_book}.

Let $N \cong \bZ^n$ be a rank $n$ lattice. We denote by $M := \text{Hom}(N,\bZ)$ the dual of $N$ and by $\langle\cdot,\cdot\rangle: M \times N \to \bZ$ the natural pairing. For any $\bZ$-module $R$, we also let $N_R := N \otimes_\bZ R$ and $M_R := M \otimes_\bZ R$.

Let $X = X_\Sigma$ be an n-dimensional toric variety defined by a fan $\Sigma$ supported in the real vector space $N_\bR = N \otimes_\bZ \bR$. The variety $X$ admits an action by the algebraic torus $T_N^\bC := N \otimes_\bZ \bC^\times \cong (\bC^\times)^n$, whence its name ``toric variety'', and it contains a Zariski open dense orbit $U_0 \subset X$ on which $T_N^\bC$ acts freely.

Denote by $v_1, \ldots, v_m \in N$ the primitive generators of rays in $\Sigma$, and by
$$D_1, \ldots, D_m \subset X$$
the corresponding toric prime divisors. We have the {\em fan sequence} associated to $X$:
\begin{equation*}
0\longrightarrow \mathbb{K}:=\text{Ker}(\phi) \longrightarrow \bZ^m \overset{\phi}{\longrightarrow} N\longrightarrow 0,
\end{equation*}
where the fan map $\phi: \bZ^d \to N$ is defined by $\phi(e_i) = v_i$, and $\{e_1,\ldots,e_m\}$ is the standard basis of $\bZ^d$. Applying $\text{Hom}(-,\bZ)$ to above gives the {\em divisor sequence}:
\begin{equation}\label{eqn:divisor_sequence}
0\longrightarrow M \longrightarrow \bigoplus_{i=1}^m \bZ D_i \longrightarrow \mathbb{K}^\vee \longrightarrow 0.
\end{equation}

In this note, we will always assume that $X = X_\Sigma$ is smooth\footnote{Most of the results described here have natural generalization to toric orbifolds, i.e. toric varieties with at most quotient singularities. We will make comments where appropriate.} and the support $|\Sigma|$ of the fan $\Sigma$ is convex and of full dimension in $N_\bR$. In this case, we have
$$\mathbb{K} \cong H_2(X;\bZ),\quad \mathbb{K}^\vee \cong \text{Pic}(X) = H^2(X;\bZ),$$
and both are of rank $r := m-n$. We let $\{p_1,\ldots,p_r\}$ be a nef basis of $H^2(X;\bZ)$ and let $\{\gamma_1,\ldots,\gamma_r\}\subset H_2(X;\bZ)$ be the dual basis.

\begin{defn}
The {\em complex moduli} of the mirror of $X = X_\Sigma$ is defined as
\begin{equation*}
\check{\mathcal{M}}_B := \mathbb{K}^\vee \otimes_{\bZ} \bC^\times.
\end{equation*}
\end{defn}
We will denote by $t = (t_1, \ldots, t_r)$ the set of coordinates on $\check{\mathcal{M}}_B \cong (\bC^\times)^r$ with respect to the nef basis $\{p_1,\ldots,p_r\} \subset H^2(X;\bZ)$, so that $t = 0$ is a {\em large complex structure limit} of the mirror complex moduli.

Suppose further that $X = X_\Sigma$ is K\"ahler. Recall that the K\"ahler cone of a K\"ahler manifold $X$ is the open cone of K\"ahler classes in $H^2(X;\bR)$.
\begin{defn}
The {\em complexified K\"ahler moduli} of $X = X_\Sigma$ is defined as
\begin{equation*}
\mathcal{M}_A := \mathcal{K}_X \oplus 2\pi\mathbf{i} \left(H^2(X;\bR) / H^2(X;\bZ)\right) \subset \mathbb{K}^\vee \otimes_\bZ \bC^\times,
\end{equation*}
where $\mathcal{K}_X \subset H^2(X;\bR)$ denotes the K\"ahler cone of $X$.
\end{defn}
We will denote by $q = (q_1, \ldots, q_r) \in \mathcal{M}_A$ the complexified K\"ahler parameters defined by
$$q_a = \exp \left( - \int_{\gamma_a} \omega_\bC \right),$$
where $\omega_\bC = \omega + 2\pi\mathbf{i} \beta \in H^2(X;\bC)$ is a complexified K\"ahler class; here the imaginary part $\beta \in H^2(X;\bR)$ is the so-called {\em $B$-field}.

A variety $X$ is called {\em Calabi-Yau} if it is Gorenstein and its canonical line bundle $K_X$ is trivial. By this definition, a toric variety $X = X_\Sigma$ is Calabi-Yau if and only if there exists a lattice point $u \in M$ such that $\langle u, v_i\rangle = 1$ for $i = 1, \ldots, m$ \cite{CLS_toric_book}. This is in turn equivalent to the existence of $u\in M$ such that the corresponding character $\chi^u \in \text{Hom}(T_M^\bC, \bC^\times)$ defines a {\em holomorphic} function on $X$ with simple zeros along each toric prime divisor $D_i$ and which is non-vanishing elsewhere. In such a case, $X$ is necessarily noncompact.

\section{Mirror symmetry for compact toric varieties}

\subsection{Landau-Ginzburg models as mirrors}

Through the works of Batyrev \cite{Batyrev93}, Givental \cite{Givental95, Givental96, Givental98}, Kontsevich \cite{Kontsevich-ENS98}, Hori-Vafa \cite{Hori-Vafa00} and many others, mirror symmetry has been extended beyond the Calabi-Yau setting. In that situation, the mirror geometry is believed to be given by a {\em Landau-Ginzburg model} which is a pair $(\check{X}, W)$ consisting of a K\"ahler manifold $\check{X}$ and a holomorphic function $W:\check{X} \to \bC$ called the {\em superpotential} of the model \cite{Vafa91, Witten93}. Note that when $W$ is nonconstant, the mirror manifold $\check{X}$ is necessarily noncompact.

Many interesting mathematical consequences can be deduced from this mirror symmetry. For instance, Kontsevich has extended his homological mirror symmetry conjecture to this setting \cite{Kontsevich-ENS98}. On the other hand, it is also predicted that genus 0 closed Gromov-Witten invariants of $X$ are encoded in the deformation or {\em unfolding} of the superpotential $W$; in particular, there should be an isomorphism of {\em Frobenius algebras}
\begin{equation}\label{eqn:QH=Jac}
QH^*(X) \cong Jac (W)
\end{equation}
between the {\em small quantum cohomology ring} $QH^*(X)$ of $X$ and the {\em Jacobian ring} $Jac(W) := \mathcal{O}(\check{X}) / J_W$ of $W$, where $J_W$ denotes the Jacobian ideal.

Let $X = X_\Sigma$ be an $n$-dimensional compact toric K\"ahler manifold defined by a fan $\Sigma$ in $N_\bR \cong \bR^n$. The mirror of $X$ is conjecturally given by a Landau-Ginzburg model $(\check{X}, W)$, where $\check{X}$ is simply the algebraic torus $T_M^\bC := M \otimes_\bZ \bC^\times \cong (\bC^\times)^n$ so that the Jacobian ring is simply
\begin{equation}\label{eqn:Jac}
Jac (W) = \frac{\bC[z_1^{\pm 1}, \ldots, z_n^{\pm 1}]}{\left\langle z_1 \frac{\partial W}{\partial z_1}, \ldots, z_n \frac{\partial W}{\partial z_n}\right\rangle}.
\end{equation}
In a more canonical way, tensoring the divisor sequence \eqref{eqn:divisor_sequence} with $\bC^\times$ gives the exact sequence
\begin{equation*}
0\longrightarrow T_M^\bC \longrightarrow (\bC^\times)^m \longrightarrow \check{\mathcal{M}}_B \longrightarrow 0,
\end{equation*}
and the mirror is given by the family of algebraic tori $(\bC^\times)^m \to \check{\mathcal{M}}_B$.

In \cite{Hori-Vafa00}, Hori and Vafa argued using physical arguments that the superpotential should be given by a family of Laurent polynomials $W^\text{HV}_t \in \bC[z_1^{\pm1},\ldots,z_n^{\pm1}]$ over the mirror complex moduli $\check{\mathcal{M}}_B$ whose Newton polytope is the convex hull of $v_1,\ldots,v_m$.
More concretely, the superpotential can be explicitly written as
\begin{equation*}
W^\text{HV}_t = \sum_{i=1}^m \check{C}_i z^{v_i},
\end{equation*}
where $z^{v}$ denotes the monomial $z_1^{v^1}\cdots z_n^{v^n}$ for $v = (v^1, \ldots, v^n) \in N \cong \bZ^n$, and the coefficients $\check{C}_i \in \bC$ are constants subject to the constraints
$$t_a = \prod_{i=1}^m \check{C}_i^{D_i\cdot \gamma_a},\quad a = 1, \ldots, r;$$
here recall that $t = (t_1, \ldots, t_r)$ are complex coordinates on $\check{\mathcal{M}}_B \cong (\bC^\times)^r$.\footnote{Scaling of the coordinates $z_1,\ldots,z_n$ gives the {\em same} superpotential, so the constants $\check{C}_1,\ldots,\check{C}_m$ are effectively determined by the parameters $t_1,\ldots,t_r$.}
We call $\left(\check{X}, W^{\text{HV}}_t\right)$ the {\em Hori-Vafa mirror} of $X$.

\begin{example}
The fan of $\bP^n$ has rays generated by the lattice points $v_i = e_i$ for $i = 1,\ldots,n$ and $v_0 = -\sum_{i=1}^n e_i$ in $N$, where $\{e_1,\ldots,e_n\}$ is the standard basis of $N = \bZ^n$. So its Hori-Vafa mirror is the family of Laurent polynomials
$$W^\text{HV}_t = z_1 + \cdots + z_n + \frac{t}{z_1\cdots z_n},$$
where $(z_1,\ldots,z_n)$ are coordinates on $\check{X} \cong (\bC^\times)^n$.
\end{example}

\begin{example}
For another example, let us consider the Hirzebruch surface $\mathbb{F}_k = \bP(\mathcal{O}_{\bP^1}(k)\oplus\mathcal{O}_{\bP^1})$ ($k \in \bZ_{\geq 0}$) whose fan has rays generated by the lattice points
$$v_1 = (1,0), v_2 = (0,1), v_3 = (-1,-k), v_4 = (0,-1) \in N = \bZ^2.$$
Its Hori-Vafa mirror is the family of Laurent polynomials
\begin{equation}\label{eqn:W^HV_Fk}
W^\text{HV}_t = z_1 + z_2 + \frac{t_1t_2^k}{z_1z_2^k} + \frac{t_2}{z_2},
\end{equation}
where $(z_1,z_2)$ are coordinates on $\check{X} \cong (\bC^\times)^2$.
\end{example}

In terms of the Hori-Vafa mirror, the isomorphism \eqref{eqn:QH=Jac} can be established as a consequence of the mirror theorem of Givental \cite{Givental98} and Lian-Liu-Yau \cite{LLY-III} in the case when $X$ is {\em semi-Fano}, i.e. when the anti-canonical divisor $-K_X$ is nef.

\subsection{The SYZ construction}

It is natural to ask whether the SYZ proposal can be generalized to non--Calabi-Yau cases as well. Auroux's pioneering work \cite{Auroux07} gave a satisfactory answer to this question. He considered pairs $(X,D)$ consisting of a compact K\"ahler manifold $X$ together with a simple normal crossing anti-canonical divisor $D \subset X$. A defining section of $D$ gives a meromorphic top form on $X$ with only simple poles along $D$ and non-vanishing elsewhere, so it defines a holomorphic volume form on the complement $X \setminus D$ (making it Calabi-Yau). Hence it makes sense to speak about special Lagrangian submanifolds in $X \setminus D$. If there is a special Lagrangian torus fibration $\rho:X\setminus D \to B$, we can then try to run the SYZ program to construct a mirror $\check{X}$.

When $X = X_\Sigma$ is a compact toric K\"ahler manifold, a canonical choice of $D$ is the union of all the toric prime divisors $D_1,\ldots,D_m$. By equipping $X$ with a toric K\"ahler structure $\omega$, we have an action on $(X,\omega)$ by the torus $T_N := N \otimes_\bZ U(1) \cong T^n$ which is Hamiltonian. The associated moment map
$$\rho: X \to \Delta,$$
where $\Delta \subset M_\bR$ denotes the moment polytope, is then a nice Lagrangian torus fibration. We view the polytope base $\Delta$ as a tropical affine manifold with boundary $\partial \Delta$. Note that the restriction of $\rho$ to the complement $U_0 = X\setminus D$ is a Lagrangian $T^n$-fibration {\em without} singular fibers. Indeed,
$$\rho|_{U_0}: U_0 \to \Delta^\circ,$$
where $\Delta^\circ$ denotes the interior of $\Delta$, is a trivial torus bundle and hence admits a section.

As in the semi-flat case, there is then a symplectomorphism
$$U_0 = X \setminus D \cong T^*\Delta^\circ/\Lambda^\vee,$$
where $\Lambda^\vee$ is the trivial lattice bundle $\Delta^\circ \times N$. The {\em SYZ mirror} $\check{X}$ is therefore simply given by taking the fiberwise dual:
$$\check{X} = T\Delta^\circ/\Lambda,$$
where $\Lambda$ is the dual lattice bundle $\Delta^\circ \times M$.

We remark that the SYZ mirror manifold $\check{X}$ is a {\em bounded domain} inside the algebraic torus $TM_\bR/ \Lambda = T_M^\bC \cong (\bC^\times)^n$ instead of the whole space \cite{Auroux07, Chan-Leung10a}. While this is natural from the point of view of SYZ, Hori and Vafa \cite{Hori-Vafa00} have proposed a renormalization procedure to enlarge the mirror and recover the whole algebraic torus; see e.g. \cite[Section 4.2]{Auroux07} for a discussion.

How about the superpotential $W$? The original SYZ proposal \cite{SYZ96} suggests that instanton corrections should stem from holomorphic disks in $X$ with boundaries lying on Lagrangian torus fibers of $\rho$. In view of this, $W$ should be closely related to counting of such holomorphic disks. This turns out to be a natural guess also from the {\em Floer-theoretic} viewpoint. Indeed, the superpotential $W$ is precisely the mirror of Fukaya-Oh-Ohta-Ono's {\em obstruction chain} $\mathfrak{m}_0$ in the Floer complex of a Lagrangian torus fiber of the moment map $\rho$ \cite{FOOO-book, FOOO-toricI, FOOO-toricII, FOOO-toricIII}, which can be expressed in terms of disk counting invariants or {\em genus 0 open Gromov-Witten invariants}, which we now review.

Let $L \subset X$ be a Lagrangian torus fiber of the moment map $\rho$. We fix a relative homotopy class $\beta \in \pi_2(X, L)$ and consider the moduli space $\overline{\mathcal{M}}_k(L;\beta)$ of stable maps from genus 0 bordered Riemann surfaces with connected boundary and $k$ boundary marked points respecting the cyclic order which represent the class $\beta$. This moduli space is a stable map compactification of the space of holomorphic disks, i.e. holomorphic embeddings $(D^2,\partial D^2) \hookrightarrow (X,L)$. We call elements in $\overline{\mathcal{M}}_k(L;\beta)$ {\em stable disks}, the domains of which are in general configurations of sphere and disk bubbles.

Fukaya, Oh, Ohta and Ono \cite{FOOO-toricI} showed that the moduli space $\overline{\mathcal{M}}_k(L;\beta)$ admits a {\em Kuranishi structure} with virtual dimension
$$ \text{vir.}\dim_\bR \overline{\mathcal{M}}_k(L;\beta) =  n+\mu(\beta)+k-3,$$
where $\mu(\beta)$ is the Maslov index of $\beta$. They also constructed, using perturbations by torus-equivariant multi-sections, a {\em virtual fundamental chain} $[\overline{\mathcal{M}}_k(L;\beta)]^\text{vir}$ \cite{FOOO-toricI}, which becomes a {\em cycle} in the case $k=1$ and $\mu(\beta) = 2$ (roughly speaking, this is because $\mu(\beta) = 2$ is the minimal Maslov index and there cannot be disk bubbles in the domains of stable disks).
\begin{defn}\label{defn:openGW}
The {\em genus 0 open Gromov-Witten invariant} $n_\beta$ is defined as
$$n_\beta := ev_*([\overline{\mathcal{M}}_1(L;\beta)]^\text{vir}) \in H_n(L;\bQ) \cong \bQ,$$
where $ev: \overline{\mathcal{M}}_1(L;\beta) \to L$ is the evaluation at the (unique) boundary marked point.
\end{defn}
It was shown in \cite{FOOO-toricI} that the numbers $n_\beta$ are independent of the choice of perturbations, so they are invariants for the pair $(X,L)$ (hence the name ``open Gromov-Witten invariants''). In general, these invariants are very difficult to compute because the moduli problem can be highly obstructed. We will come back to this in the next section.

In terms of the open Gromov-Witten invariants $n_\beta$, we have the following
\begin{defn}\label{defn:Lagr_Floer_W}
The {\em Lagrangian Floer superpotential} is defined as
\begin{equation}\label{eqn:W^LF}
W^{\text{LF}}_q := \sum_{\substack{\beta \in \pi_2 (X, L)\\ \mu(\beta) = 2}} n_\beta Z_{\beta}.
\end{equation}
Here, we regard $\check{X}$ as the moduli space of pairs $(L,\nabla)$ where $L$ is a Lagrangian torus fiber of the moment map $\rho$ and $\nabla$ is a flat $U(1)$-connection over $L$, and $Z_\beta$ is the function (in fact a monomial) explicitly defined by
$$Z_{\beta}(L,\nabla) = \exp\left( -\int_\beta \omega_q \right)\text{hol}_{\nabla}(\partial \beta ),$$
where $\text{hol}_\nabla(\partial \beta)$ denotes the holonomy of the connection $\nabla$ around the loop $\partial\beta$ and $q = (q_1, \ldots, q_r) \in \mathcal{M}_A$ are complexified K\"ahler parameters.
\end{defn}

\begin{remark}
In general, the sum \eqref{eqn:W^LF} involves infinitely many terms. In order to make sense of it, one needs to regard the $q_a$'s as formal Novikov variables and regard \eqref{eqn:W^LF} as a formal power series over the Novikov ring. However, for simplicity of exposition, we will always assume that the power series in this paper all converge so that we can work over $\bC$. (We will see later that these series indeed converge in the semi-Fano case.) Then we can regard $W^\text{LF}_q$ as a family of holomorphic functions parameterized by $q = (q_1, \ldots, q_r) \in \mathcal{M}_A$.
\end{remark}

In \cite{Chan-Leung10a}, mirror symmetry for toric Fano manifolds was used as a testing ground to see how useful fiberwise Fourier--type transforms, or SYZ transforms, could be in the investigation of the geometry of mirror symmetry.

Before we proceed, let us look at Fourier transforms on a single torus $T=V/\Lambda$ more closely.
Recall that the semi-flat SYZ transform $\mathcal{F}^\text{semi-flat}: \Omega^{0}\left(T\right)\rightarrow\Omega^{n}\left(T^*\right)$ is given by $$\mathcal{F}^\text{semi-flat}\left(\phi\right)=\int_{T}\phi\left(y^{1},\cdots,y^{n}\right)e^{\mathbf{i}\sum_{j}dy_{j}\wedge dy^{j}},$$
while the usual Fourier transform (or series) gives a function $\hat{\phi}:\Lambda^\vee \rightarrow\mathbb{C}$ on the dual lattice $\Lambda^\vee$ by
$$
\hat{\phi}\left(m_{1},...,m_{n}\right) = \int_{T}\phi\left(y_{1},\cdots,y_{n}\right)  e^{\mathbf{i}\sum_{j}m_{j}y^{j}} dy_{1}\wedge\cdots\wedge dy_{n}.
$$
It is natural to combine these two transforms to define a transformation which transforms differential forms on $T \times \Lambda$ to those on $T^* \times \Lambda^\vee$.

Note that $T \times \Lambda$ is the space of geodesic (or affine) loops in $T$ sitting inside the loop space of $T$, i.e. $$T\times\Lambda=\mathcal{L}_{\text{geo}}T\subset\mathcal{L}T,$$
and the loop space certainly plays an important role in string theory.
We are going to describe how such a transformation $\mathcal{F}^\text{SYZ}$, called the \textit{SYZ transform}, interchanges symplectic and complex geometries, with quantum corrections included, when $X$ is a compact toric Fano manifold.

Let us consider the open dense subset $U_0 := X\setminus D$, which is the union of all the Lagrangian torus fibers of the moment map.
Recall that we have an identification $U_0 \cong T^{\ast}\Delta^\circ/\Lambda^\vee$ as symplectic manifolds. Next we consider the space
$$\tilde{X} := U_0 \times \Lambda \subset \mathcal{L} X$$
of fiberwise geodesic/affine loops in $U_0$. On $\tilde{X}$, we have the {\em quantum-corrected} symplectic structure
$$\tilde{\omega} = \omega + \Phi,$$
where the values of $\Phi$ are the Fourier modes of the Lagrangian Floer superpotential $W^{\text{LF}}_q$, so they are generating functions of the genus 0 open Gromov-Witten invariants $n_\beta$ in Definition \ref{defn:openGW} which count (virtually) holomorphic disks bounded by the Lagrangian torus fibers of the moment map $\rho$.

An explicit SYZ transform $\mathcal{F}^\text{SYZ}$ can then be constructed by {\em combining} the semi-flat SYZ transform $\mathcal{F}^\text{semi-flat}$ with fiberwise Fourier series, which can be applied to understand the geometry of mirror symmetry:
\begin{thm}[\cite{Chan-Leung10a}, Theorem 1.1]
\hfill
\begin{itemize}
\item[(1)]
The SYZ transforms interchange the quantum-corrected symplectic structure on $\tilde{X}$ with the holomorphic volume form of the pair $(\check{X},W)$, i.e. we have
\begin{align*}
\mathcal{F}^\text{SYZ}\left(e^{\mathbf{i}(\omega + \Phi)}\right) & = e^{W}\check{\Omega},\\
\left(\mathcal{F}^\text{SYZ}\right)^{-1}\left(e^{W}\check{\Omega}\right) & = e^{\mathbf{i}(\omega + \Phi)}.
\end{align*}
\item[(2)]
The SYZ transform induces an isomorphism
$$\mathcal{F}^\text{SYZ}: QH^*\left( X, \omega_q \right) \overset{\cong}{\rightarrow} Jac\left( W^{\text{LF}}_q \right)$$
between the small quantum cohomology ring of $X$ and the Jacobian ring of the Lagrangian Floer superpotential $W^{\text{LF}}_q$.\footnote{Note that here $W^{\text{LF}}_q$ is defined only on a bounded domain inside $T_M^\bC \cong (\bC^\times)^n$ but we still use the formula \eqref{eqn:Jac} to define the Jacobian ring $Jac\left( W^{\text{LF}}_q \right)$.}
\end{itemize}
\end{thm}

\begin{figure}[ht]
\setlength{\unitlength}{1mm}
\begin{picture}(100,35)
\curve(25,17, 43,17) \curve(25,17, 25,35) \curve(25,17, 10,2) \put(24.1,16.1){$\bullet$}
\put(52,16){glued from} \curve(89,17, 108,17) \put(88,16.1){$\bullet$} \curve(87,19, 87,37) \put(86.1,18){$\bullet$} \curve(86,16, 71,1) \put(85,15){$\bullet$}
\end{picture}\caption{A tropical curve as a gluing of tropical disks.}\label{fig:tropical}
\end{figure}
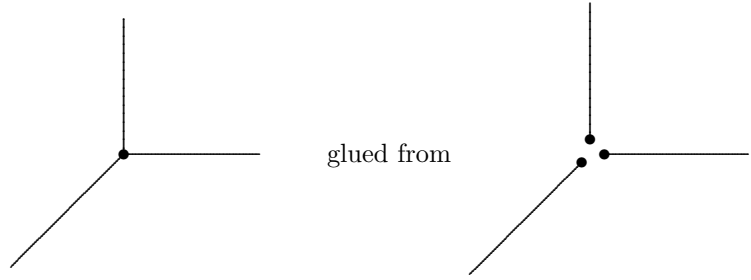

The second part of this theorem was proved using the idea that holomorphic spheres can be described as suitable gluing of holomorphic disks. More precisely, we pass to the tropical limit (see Mikhalkin \cite{Mikhalkin06} for a nice introduction to tropical geometry) and observe that a tropical curve can be obtained as a gluing of tropical disks. For example, a line $\bP^1$ in $\bP^2$ is obtained as a gluing of three disks, and this can be made precise in the tropical limit (see Figure \ref{fig:tropical}). This observation was later generalized and used by Gross \cite{Gross10} in his study of mirror symmetry for the big quantum cohomology of $\bP^2$ via tropical geometry; see also \cite{Gross_book}.

The first part of the above theorem can in fact be generalized to {\em any} compact toric manifold by appropriately defining the function $\Phi$ as a power series over certain Novikov rings; the second part was proved for semi-Fano toric surfaces by direct computations in \cite{Chan-Lau10}, and in general it follows from the main theorem in Fukaya-Oh-Ohta-Ono \cite{FOOO-toricIII}.




\subsection{Open mirror theorems}\label{sec:open_mirror_thm}

For semi-Fano toric manifolds, the mirror theorems of Givental \cite{Givental98} and Lian-Liu-Yau \cite{LLY-III} imply that there is an isomorphism
\begin{equation*}
QH^*(X,\omega_q) \cong Jac\left( W^{\text{HV}}_{t(q)} \right),
\end{equation*}
where $t(q) = \psi^{-1}(q)$ is the inverse of the {\em toric mirror map}
\begin{align*}
\psi: \check{\mathcal{M}}_B & \to \mathcal{M}_A,\\
t=(t_1, \ldots, t_r) & \mapsto \psi(t) = (q_1(t), \ldots, q_r(t)),
\end{align*}
which, by definition, is given by the $1/z$-coefficient of the combinatorially defined, cohomology-valued {\em $I$-function}:
\begin{equation*}
I_X(t,z) = e^{\sum_{a=1}^r p_a \log t_a/z} \sum_{d\in H_2^\text{eff}(X;\bZ)} t^d \prod_{i=1}^m \frac{\prod_{k = -\infty}^0 (D_i + kz)}{\prod_{k = -\infty}^{D_i\cdot d}(D_i + kz)},
\end{equation*}
where $t^d = \prod_{a=1}^{r} t_a^{p_a\cdot d}$ and $D_i$ is identified with its cohomology class in $H^2(X)$.
The toric mirror map $\psi$, which gives a local isomorphism between the mirror complex moduli $\check{\mathcal{M}}_B$ and the complexified K\"ahler moduli $\mathcal{M}_A$ of $X$ near the large complex structure and volume limits respectively, can also be obtained by solving a system of PDEs known as {\em Picard-Fuchs equations}.

On the other hand, instanton corrections can be realized using the genus 0 open Gromov-Witten invariants $n_\beta$. In terms of the Lagrangian Floer superpotential (Definition \ref{defn:Lagr_Floer_W})
$$W^{\text{LF}}_q = \sum_{\substack{\beta \in \pi_2 (X, L)\\ \mu(\beta) = 2}} n_\beta Z_{\beta},$$
Fukaya, Oh, Ohta and Ono proved in \cite{FOOO-toricIII} that there is an isomorphism of Frobenius algebras\footnote{In fact what Fukaya, Oh, Ohta and Ono obtained in \cite{FOOO-toricIII} is a ring isomorphism between the {\em big} quantum cohomology ring of $X$ and the Jacobian ring of the so-called {\em bulk-deformed} superpotential \cite{FOOO-toricII}. For the purpose of this note, we restrict ourselves to the small case.}
\begin{equation*}
QH^*(X,\omega_q) \cong Jac\left( W^{\text{LF}}_q \right).
\end{equation*}

While the superpotentials $W^{\text{HV}}$ and $W^{\text{LF}}$ originate from different approaches, they give rise to essentially the same mirror symmetry statement. This leads us to the following conjecture:
\begin{conj}\label{conj:W_HV=W_LF}
Let $X$ be a semi-Fano toric manifold, i.e. its anti-canonical divisor $-K_X$ is nef. Let $W^{\text{HV}}$ and $W^{\text{LF}}$ be the Hori-Vafa and Lagrangian Floer superpotentials of $X$ respectively. Then we have the equality
\begin{equation}
W^{\text{HV}}_{t(q)} = W^{\text{LF}}_q
\end{equation}
via the inverse mirror map $t(q) = \psi^{-1}(q)$.
\end{conj}

In \cite{Cho-Oh06}, Cho and Oh classified all embedded holomorphic disks in $X$ bounded by a fixed Lagrangian torus fiber $L$ of the moment map. Applying this classification, they computed all the invariants $n_\beta$ in the Fano case:
$$ n_\beta = \left\{
\begin{array}{ll}
1 & \text{if $\beta = \beta_i$ for some $i=1,\ldots,m$},\\
0 & \text{otherwise}.
\end{array}\right.$$
Here, the $\beta_i$'s are the so-called {\em basic disk classes}: for each $i \in \{1, \ldots, m\}$, the relative homotopy class $\beta_i$ is represented by an embedded Maslov index 2 holomorphic disk in $X$ with boundary on $L$ which intersects the toric divisor $D_i$ at exactly one point in its interior; it can be shown that $\{\beta_1, \ldots ,\beta_m\}$ form a $\bZ$-basis of the relative homotopy group $\pi_2(X,L)$ \cite{Cho-Oh06}.

The results of Cho and Oh give an explicit formula for $W^\text{LF}$ and also a direct verification of the above conjecture in the Fano case. Note that in this case the mirror map is trivial.
In fact, Cho and Oh proved that $n_{\beta_i} = 1$ for all $i$ for {\em any} compact toric manifold. However, for a general homotopy class $\beta \in \pi_2(X,L)$ which maybe represented by a stable map with disk and/or sphere bubbles, the invariants $n_\beta$ are usually very difficult to compute.

Nevertheless, in the semi-Fano case, since the first Chern number of any effective curve is nonnegative, only sphere bubbles can appear in Maslov index 2 stable disks. Therefore, $n_\beta$ is nonzero only when $\beta$ is of the form
$$\beta = \beta_i + d,$$
where $\beta_i$ is a basic disk class and $d \in H_2^\text{eff}(X;\bZ)$ is an effective curve class with first Chern number
$$c_1(d) := c_1(X) \cdot d = 0.$$
In particular, the Lagrangian Floer superpotential is a Laurent polynomial which can be expressed nicely as:
\begin{equation*}
W^\text{LF}_q = \sum_{i=1}^m C_i(1+\delta_i(q)) z^{v_i},
\end{equation*}
where
$$ 1 + \delta_i(q) := \sum_{\substack{\alpha\in H_2^\text{eff}(X;\bZ);\\ c_1(\alpha) = 0}} n_{\beta_i + \alpha} q^\alpha$$
is a generating function of genus 0 open Gromov-Witten invariants, and the coefficients $C_i \in \bC$ are constants subject to the constraints
$$q_a = \prod_{i=1}^m C_i^{D_i\cdot \gamma_a}, \quad a = 1, \ldots, r,$$
where $q = (q_1, \ldots, q_r) \in \mathcal{M}_A$ are the complexified K\"ahler parameters.

This suggests an equivalent formulation of Conjecture \ref{conj:W_HV=W_LF} in terms of the so-called {\em SYZ map}:
\begin{defn}
The {\em SYZ map} is the map
\begin{align*}
\phi: \mathcal{M}_A & \to \check{\mathcal{M}}_B,\\
q = (q_1, \ldots, q_r) & \mapsto \phi(q) = (t_1(q), \ldots, t_r(q))
\end{align*}
defined by
\begin{equation*}
t^a(q) = q_a \prod_{i=1}^m \left( 1+\delta_i(q) \right)^{D_i\cdot \gamma_a}
\end{equation*}
for $a = 1,\ldots, r$.
\end{defn}

Then Conjecture \ref{conj:W_HV=W_LF} is equivalent to the following
\begin{conj}\label{conj:SYZ_map=mirror_map}
The SYZ map $\phi$ is inverse to the mirror map $\psi$.
\end{conj}
In particular this gives an {\em enumerative meaning} to inverses of mirror maps, and explains the mysterious integrality property of mirror maps observed earlier in \cite{Lian-Yau96, Lian-Yau98, Krattenthaler-Rivoal10, Zhou12}

To prove the above conjectures, we need to compute the invariants $n_\beta$. The first known non-Fano example is the Hirzebruch surface $\mathbb{F}_{2}$, whose Lagrangian Floer superpotential was first computed by Auroux \cite{Auroux09} using wall-crossing and degeneration techniques (which work for $\mathbb{F}_{3}$ as well) and then by Fukaya-Oh-Ohta-Ono in \cite{FOOO12} by employing their big machinery \cite{FOOO-book}. The result is given by the Laurent polynomial
\begin{equation}\label{eqn:W^LF_F2}
W^\text{LF}_q = z_1 + z_2 + \frac{q_1 q_2^2}{z_1 z_2^2} + \frac{q_2(1+q_1)}{z_2}.
\end{equation}
Comparing \eqref{eqn:W^HV_Fk} (for $k=2$) and \eqref{eqn:W^LF_F2}, we see that, besides the basic disks, there is one extra stable disk whose domain contains a sphere bubble which is a copy of the $\left(-2\right)$-curve in $\mathbb{F}_{2}$; see Figure \ref{bubble_F2} below.

\begin{figure}[htp]
\begin{center}
\includegraphics{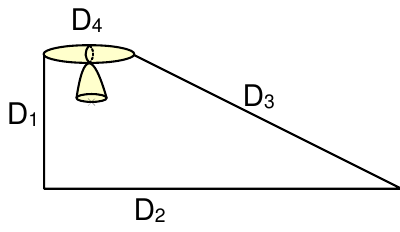}
\end{center}
\caption{The stable disk in $\mathbb{F}_{2}$ as a union of a disk and the sphere $D_4$.} \label{bubble_F2}
\end{figure}

In \cite{Chan11}, the invariants $n_\beta$ were computed for any toric $\bP^1$-bundle of the form $X = \bP(K_Y \oplus \mathcal{O}_Y)$, where $Y$ is a compact toric Fano manifold and $K_Y$ denotes (by abuse of notations) the canonical line bundle over $Y$; note that this class of examples includes the Hirzebruch surface $\mathbb{F}_2 = \bP(K_{\bP^1} \oplus \mathcal{O}_{\bP^1})$. The computation is done by establishing the following equality \cite[Theorem 1.1]{Chan11}:
\begin{equation}\label{eqn:open_closed}
n_\beta = \langle [\text{pt}] \rangle^X_{0,1,\bar{\beta}},
\end{equation}
between open and closed Gromov-Witten invariants. The right-hand side of the equality is the genus 0 closed Gromov-Witten invariant $\langle [\text{pt}] \rangle^X_{0,1,\bar{\beta}}$, where $[\text{pt}] \in H^{2n}(X;\bQ)$ is the class of a point and $\bar{\beta}$ is the class of the closed stable map obtained by capping off the disk component of a stable disk representing $\beta$. This equality, together with the observation that the invariants $\langle [\text{pt}] \rangle^X_{0,1,\bar{\beta}}$ appear in a certain coefficient of Givental's {\em $J$-function} and also a computation of the mirror map \cite{CLT11}, prove both Conjectures \ref{conj:W_HV=W_LF} and \ref{conj:SYZ_map=mirror_map} for the class of examples $X=\bP(K_Y \oplus \mathcal{O}_Y)$.

In \cite{Chan-Lau10}, the open Gromov-Witten invariants $n_\beta$ and hence the Lagrangian Floer superpotential $W^\text{LF}$ were computed for all semi-Fano toric surfaces. For these surfaces, every toric divisor (curve) has self-intersection number at least $-2$ (indeed, this is equivalent to the surface being semi-Fano). Stable disks with nontrivial sphere bubbles can only appear along a chain of $\left(-2\right)$-toric curves. The above open/closed equality then allows us to identify the invariant $n_\beta$ with a certain closed {\em local} Gromov-Witten invariant, which was previously computed by Bryan and Leung \cite{Bryan-Leung00}. We thus obtain explicit formulas for both $n_\beta$ and $W^\text{LF}$; see \cite[Theorem 1.2 and Appendix A]{Chan-Lau10}. As an application, an explicit presentation for the small quantum cohomology ring of such a surface was also obtained; see \cite[Section 4]{Chan-Lau10} for details.

As an example, consider the semi-Fano toric surface $X = X_\Sigma$ defined by the fan and polytope shown in Figure \ref{semi-fano_surf} below.

\begin{figure}[ht]
\begin{center}
\includegraphics{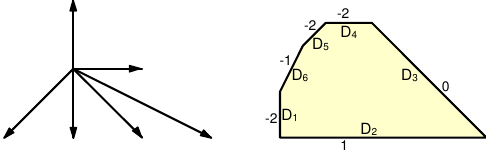}
\end{center}
\caption{The fan $\Sigma$ (left) and the polytope $\Delta$ (right) defining $X$; the number besides each edge is the self-intersection number of the corresponding toric divisor.}\label{semi-fano_surf}
\end{figure}

\noindent Then the superpotential $W^\text{LF}$ is given by
\begin{eqnarray*}
W^\text{LF} = & (1 + q_1)z_1 + z_2 + \frac{q_1q_2q_3^2q_4^3}{z_1z_2} + (1+q_2+q_2q_3)\frac{q_1q_3q_4^2}{z_2}\\
              & \qquad + (1+q_3+q_2q_3)\frac{q_1q_4z_1}{z_2} + \frac{q_1z_1^2}{z_2},
\end{eqnarray*}
where $q_l=\exp(-t_l)$, $l=1,\ldots,4$, are complexified K\"ahler parameters on $\mathcal{M}_A$.

For general semi-Fano toric manifolds, the major difficulty is to find a suitable space to replace $X$ on the right-hand side of the open/closed equality \eqref{eqn:open_closed} so that it remains valid. In \cite{CLLT12}, it was discovered that a natural choice was given by the so-called {\em Seidel space} $E_i$, an $X$-bundle over $\bP^1$ corresponding to the $\bC^\times$-action defined by the divisor $D_i$ which has been used by Seidel \cite{Seidel97} to construct the {\em Seidel representation} -- a nontrivial action of $\pi_1(\text{Ham}(X,\omega))$ on $QH^*(X,\omega)$.

By establishing an open/closed equality similar to \eqref{eqn:open_closed} but with $X$ replaced by $E_i$ on its right-hand side, and applying various techniques such as certain degeneration formulas, the following theorem was proved in \cite{CLLT12}.

\begin{thm}[\cite{CLLT12}, Theorem 1.2 and Corollary 6.6]
Both Conjectures \ref{conj:W_HV=W_LF} and \ref{conj:SYZ_map=mirror_map} are true for any compact semi-Fano toric manifold.
\end{thm}

We call this the {\em open mirror theorem}. When combined with the closed mirror theorem of Givental \cite{Givental98} and Lian-Liu-Yau \cite{LLY-III}, this gives the isomorphism
\begin{equation*}
QH^*(X,\omega_q) \cong Jac(W^{\text{LF}}_q),
\end{equation*}
which recovers the aforementioned result of Fukaya, Oh, Ohta and Ono \cite{FOOO-toricIII}.

More importantly, the open mirror theorem gives us explicit formulas which can effectively compute {\em all} the genus 0 open Gromov-Witten invariants $n_\beta$ in the semi-Fano case:
\begin{thm}[\cite{CLLT12}, Theorem 1.1]\label{thm:open_mirror_thm}
Let $X$ be a compact semi-Fano toric manifold. Then we have the following formula for the generating function of genus 0 open Gromov-Witten invariants:
$$1+\delta_i(q)= \exp g_i(t(q)),$$
where $t(q) = \psi^{-1}(q)$ is the inverse of the mirror map $q=\psi(t)$ and $g_i$ is a hypergeometric function explicitly defined by
\begin{equation*}
g_i(t):=\sum_{d}\frac{(-1)^{(D_i\cdot d)}(-(D_i\cdot d)-1)!}{\prod_{p\neq i} (D_p\cdot d)!} t^d;
\end{equation*}
here the summation is over all effective curve classes $d\in H_2^\text{eff}(X;\bZ)$ satisfying
$$-K_X\cdot d=0, D_i\cdot d<0 \text{ and } D_p\cdot d \geq 0 \text{ for all } p\neq i.$$
\end{thm}

As another immediate corollary, we also have the following convergence result.
\begin{thm}[\cite{CLLT12}, Theorem 1.3] \label{thm conv}
For a compact semi-Fano toric manifold $X$, the generating functions $1+\delta_i(q)$ are convergent power series, and hence the Lagrangian Floer superpotential $W^{\text{LF}}_q$ is a family of holomorphic functions parametrized by the complexified K\"ahler parameters $q = (q_1, \ldots, q_r) \in \mathcal{M}_A$.
\end{thm}

We refer the reader to the paper \cite{CLLT12} for more details. We remark that an open mirror theorem in the orbifold setting was also proved in \cite{CCLT12} for weighted projective spaces of the form $\bP(1,\ldots,1,n)$ using similar techniques; in fact, Theorem \ref{thm:open_mirror_thm} can be extended to all semi-Fano toric {\em orbifolds} as well \cite{CCLLT15}.

\section{Mirror symmetry for toric Calabi-Yau varieties}

\subsection{Physicists' mirrors}

Let $X = X_\Sigma$ be a toric Calabi-Yau variety of complex dimension $n$.
By choosing a suitable basis of $N = \bZ^n$, we may write
$$v_i = (w_i, 1) \in N = \bZ^{n-1} \oplus \bZ,$$
where $w_i \in \bZ^{n-1}$. Recall that we have assumed that $X$ is smooth and the fan $\Sigma$ has convex support. In this case, the Calabi-Yau manifold $X$ is a crepant resolution of an affine toric variety (defined by the cone $|\Sigma|$) with Gorenstein canonical singularities, which is equivalent to saying that $X$ is {\em semi-projective} \cite[p.332]{CLS_toric_book}.

An important class of examples of toric Calabi-Yau manifolds is given by total spaces of the canonical line bundles $K_Y$ over compact toric manifolds $Y$.
For example, the total space of $K_{\bP^2}=\mathcal{O}_{\bP^2}(-3)$ is a toric Calabi-Yau 3-fold whose fan $\Sigma$ has rays generated by the lattice points
$$v_0=(0,0,1), v_1=(1,0,1), v_2=(0,1,1), v_3=(-1,-1,1)\in N = \bZ^3.$$

Mirror symmetry in this setting is called {\em local mirror symmetry} because it originated from an application of mirror symmetry techniques to Fano surfaces (e.g $\bP^2$) contained inside compact Calabi-Yau manifolds and could be derived using physical arguments from mirror symmetry for compact Calabi-Yau hypersurfaces in toric varieties by taking certain limits in the complexified K\"ahler and complex moduli spaces \cite{KKV97}.

Given a toric Calabi-Yau manifold $X$ as above, its mirror is predicted to be a family of affine hypersurfaces in $\bC^2 \times (\bC^\times)^{n-1}$ \cite{Leung-Vafa98, CKYZ99, HIV00}
\begin{equation}\label{eqn:toricCY_mirror}
\check{X}_t = \left\{ (u, v, z_1, \ldots, z_{n-1}) \in \bC^2 \times (\bC^\times)^{n-1} : uv = \sum_{i=1}^{m} \check{C}_i z^{w_i} \right\},
\end{equation}
where the coefficients $\check{C}_i \in \bC$ are constants subject to the constraints
$$t_a = \prod_{i=1}^m \check{C}_i^{D_i\cdot \gamma_a},\quad a=1,\ldots,r;$$
here, again, $t = (t_1, \ldots, t_r)$ are coordinates on the mirror complex moduli $\check{\mathcal{M}}_B := \mathbb{K}^\vee \otimes_{\bZ} \bC^\times \cong (\bC^\times)^r$. $\check{X}_t$ are noncompact Calabi-Yau manifolds since
$$\check{\Omega}_t := \text{Res}\left( \frac{du \wedge dv \wedge d\log z_1 \wedge \cdots \wedge d\log z_{n-1}}{uv - \sum_{i=1}^{m} \check{C}_i z^{w_i}} \right)$$
is a nowhere vanishing holomorphic volume form on $\check{X}_t$.

\begin{example}
When $X = K_{\bP^2}$, its mirror is predicted to be
\begin{align*}
\check{X}_t = \left\{ (u,v,z_1,z_2) \in \bC^2 \times (\bC^\times)^2 : uv = 1 + z_1 + z_2 + \frac{t}{z_1z_2} \right\},
\end{align*}
where $t$ is a coordinate on the mirror complex moduli $\check{\mathcal{M}}_B \cong \bC^\times$.
\end{example}

This mirror symmetry has been a rich source of interesting examples and has numerous applications. This has resulted in a flurry of research works by both physicists and mathematicians \cite{Leung-Vafa98, CKYZ99, HIV00, Gross01, Gross-Inventiones, Takahashi01, Klemm-Zaslow01, Graber-Zaslow02, Hosono00, Hosono06, Forbes-Jinzenji05, Forbes-Jinzenji06, Konishi-Minabe10, Seidel10}.\footnote{This list is certainly not meant to be exhaustive.} In particular, just like the usual case, it has been applied to make powerful enumerative predictions.

For instance, in the 3-fold case, mirror symmetry was employed to compute the {\em local Gromov-Witten invariants}
\begin{align*}
N_{g,d}(X) := \int_{[\overline{\mathcal{M}}_{g,0}(X,d)]^\text{vir}} \textbf{1},
\end{align*}
where $\overline{\mathcal{M}}_{g,0}(X,d)$ is the moduli space of genus $g$, 0-pointed closed stable maps into $X$ with class $d \in H_2^{\text{eff}}(X;\bZ)$. These invariants are in general hard to compute due to nontrivial obstructions, but mirror symmetry can be used to compute (at least) the genus 0 invariants all at once.

To carry out the computations, one needs a {\em mirror map}
$$\psi: \check{\mathcal{M}}_B \to \mathcal{M}_A,$$
which is a local isomorphism from the complex moduli $\check{\mathcal{M}}_B$ of the mirror $\check{X}$ to the complexified K\"ahler moduli
$\mathcal{M}_A$ of $X$. Traditionally, mirror maps are defined as quotients of {\em period integrals}, which are integrals of the form
\begin{align*}
\int_\Gamma \check{\Omega},
\end{align*}
where $\Gamma \in H_3(\check{X};\bZ)$ is a middle-dimensional integral cycle.
These integrals satisfy a system of {\em Picard-Fuchs} PDEs \cite{Batyrev93}, whose solutions are explicit hypergeometric--type functions, and indeed this is usually how mirror maps are computed.

For the example of $X = K_{\bP^2}$, the Picard-Fuchs system is just a single ODE
\begin{align*}
\left[\Theta^3 + 3t\Theta(3\Theta + 1)(3\Theta + 2) \right]\Phi = 0,
\end{align*}
where $\Theta = t d/dt$.
A basis of solutions to this equation is explicitly given by
\begin{align*}
\Phi_0(t) & = 1,\\
\Phi_1(t) & = \log t + \sum_{k=1}^\infty \frac{(-1)^k}{k}\frac{(3k)!}{(k!)^3} t^k, \\
\Phi_2(t) & = (\log t)^2 + \cdots.
\end{align*}
Note that constants are the only holomorphic solutions of the Picards-Fuchs equation, and this is true for general toric Calabi-Yau manifolds.

The mirror map is then given by
$$ q = \psi(t) := \exp \left( \Phi_1(t)/\Phi_0(t) \right) = t(1 + \cdots),$$
and genus 0 local Gromov-Witten invariants $N_{0,d}(X)$ are predicted to appear as coefficients of the Taylor series expansion of $\Phi_2(\psi^{-1}(q))$ around $q = 0$.

In \cite{Graber-Zaslow02}, it was shown that, via the inverse mirror map $q \mapsto \psi^{-1}(q)$ and a coordinate change, the mirror of $K_{\bP^2}$ can be rewritten as the family
\begin{align*}
\check{X}_q = \left\{ (u,v,z_1,z_2) \in \bC^2 \times (\bC^\times)^2 : uv = f(q) + z_1 + z_2 + \frac{q}{z_1z_2} \right\},
\end{align*}
where
\begin{equation}\label{eqn:KP2_integral_series}
f(q) = \left(\psi(t)/t \right)^{\frac{1}{3}} = 1 - 2q + 5q^2 -32q^3 + 286q^4 - 3038q^5 + \cdots.
\end{equation}
You may wonder why the coefficients of $f(q)$ are all integers and whether they are counting something. We are going to see that these questions can all be answered by applying the SYZ construction, and that the coefficients of $f(q)$ are (virtual) counts of holomorphic disks.

\subsection{The SYZ construction}

To carry out the SYZ mirror construction for a toric Calabi-Yau manifold (or more generally an orbifold) $X$, we need a special Lagrangian torus fibration on $X$. This is provided by a construction due independently to Goldstein \cite{Goldstein01} and Gross \cite{Gross01}.

To begin with, recall that the lattice point $u \in M$, which defines the hyperplane containing all the generators $v_i$'s, corresponds to a holomorphic function $\chi^u: X \to \bC$ with simple zeros along each toric prime divisor $D_i \subset X$. We equip $X$ with a toric K\"ahler structure $\omega$ and consider the action by the subtorus $T^{n-1} \subset T_N \cong T^n$ which preserves $\chi^u$, or equivalently, the subtorus whose action preserves the canonical holomorphic volume form $\Omega$ on $X$. Let $\rho_0: X \to \bR^{n-1}$ be the corresponding moment map which is given by composing the $T_N$-moment map with the projection along the ray in $M_\bR$ spanned by $u$.
\begin{prop}[Goldstein \cite{Goldstein01}, Gross \cite{Gross01}]
For any nonzero constant $\epsilon \in \bC^\times$, the map defined by
\begin{align*}
\rho := \left( \rho_0, |\chi^u - \epsilon| \right): X \to B:= \bR^{n-1} \times \bR_{\geq0},
\end{align*}
is a special Lagrangian torus fibration, where the fibers are special with respect to the meromorphic volume form
$$\Omega_\epsilon := \frac{\Omega}{\chi^u - \epsilon}.$$
\end{prop}

We call $\rho$ the {\em Gross fibration}, which is {\em non-toric} in the sense that it is not the $T_N$-moment map. Its discriminant locus is easy to describe, namely, a fiber of $\rho$ is non-regular if and only if either
\begin{itemize}
\item
when it intersects nontrivially with (and hence is contained inside) the hypersurface $D_\epsilon \subset X$ defined by $\chi^u = \epsilon$, in which case the fiber is mapped to a point on the boundary $\partial B = \bR^{n-1} \times \{0\}$; or
\item
when it contains a point where the $T^{n-1}$-action is not free, i.e. when at least two of the homogeneous coordinates on $X$ vanish, in which case the fiber is mapped to the image $\Gamma$ of the codimension 2 subvariety
$$\bigcup_{i\neq j} D_i \cap D_j$$
under $\rho$.
\end{itemize}
We regard $B$ as a tropical affine manifold with boundary $\partial B$ and singularities $\Gamma$. Note that $\Gamma$ has real codimension 2 in $B$.

\begin{figure}
\begin{center}
\begin{tikzpicture}
\draw[dashed] (2,1.5) -- (7,0.75);
\draw[dashed] (2,1.5) -- (5,2.5);
\draw[red] (6.2,1.5) node[left] {$\Gamma$};
\draw[thick] (2.05,1.5) -- (1.95,1.5) node[left] {$|\epsilon|$};
\draw (1.95,0.75) node[left] {$B_-$};
\draw (1.95,2.25) node[left] {$B_+$};
\draw[->] (1.5,0.075) -- (7,-0.75);
\draw[->] (1.55,-0.15) -- (5,1);
\draw (5.2,0.2) node[left] {$\partial B$};
\draw[->] (2,-0.5) -- (2,3);
\draw [thick, red] (3,1.6) -- (4.5,1.75);
\draw [thick, red] (4.5,1.75) -- (5.15,1.9665);
\draw [thick, red] (4.5,1.75) -- (5.25,1.6375);
\draw [thick, red] (5.15,1.9665) -- (5.25,1.6375);
\draw [thick, red] (5.15,1.9665) -- (5.45,2.35);
\draw [thick, red] (5.25,1.6375) -- (6.15,1.1);
\end{tikzpicture}
\caption{The base of the Gross fibration on $K_{\bP^2}$, which is an upper half space in $\bR^3$.}\label{KP2_base}
\end{center}
\end{figure}
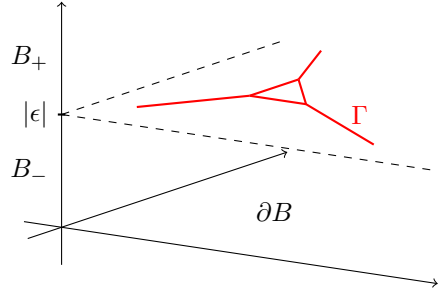

For example, the base $B$ of the Gross fibration on $X = K_{\bP^2}$ is an upper half space in $\bR^3$, and the discriminant locus is a graph $\Gamma$ which is contained in a hyperplane $H$ parallel to the boundary $\partial B$, as shown in Figure \ref{KP2_base}.

Starting with a (special) Lagrangian torus fibration $\rho: X \to B$ on a Calabi-Yau manifold, the SYZ mirror construction proceeds in several steps:
\begin{itemize}

\item[Step 1]
Over the smooth locus $B_0 := B \setminus \left( \partial B \cup \Gamma \right)$, the pre-image $X_0 := \rho^{-1}(B_0)$ can be identified with the quotient $T^*B_0/\Lambda^\vee$ by Duistermaat's action-angle coordinates.

\item[Step 2]
Define the {\em semi-flat} mirror $\check{X}_0$ as $TB_0/\Lambda$, which is not quite the correct mirror because the complex structure on $\check{X}_0$ {\em cannot} be extended further to {\em any} (partial) compactification as the monodromy of the tropical affine structure around the discriminant locus $\Gamma$ is nontrivial.

\item[Step 3]
To obtain the correct and (partially) compactified mirror $\check{X}$, we need to modify the complex structure on $\check{X}_0$ by instanton corrections, which should come from holomorphic disks in $X$ bounded by Lagrangian torus fibers of the fibration $\rho$.
\end{itemize}

Such a procedure was pioneered and first put into practice by Auroux in \cite{Auroux07, Auroux09}. Among the many interesting and motivating examples studied by him there was $X = \bC^n$, which can be regarded as the simplest example of a toric Calabi-Yau manifold. The work \cite{CLL12} gives a generalization of Auroux's construction to all toric Calabi-Yau manifolds, which we now review.

First of all, by definition, the {\em wall(s)} inside the base of a Lagrangian torus fibration is the loci of regular fibers which bound nonconstant {\em Maslov index 0} holomorphic or stable disks in $X$. For the Gross fibration on a toric Calabi-Yau manifold, there is only one wall given by the hyperplane
$$H := \bR^{n-1} \times \{|\epsilon|\} \subset B$$
parallel to the boundary $\partial B$. The wall $H$ contains the discriminant locus $\Gamma$ as a {\em tropical hypersurface}, and it divides the base $B$ into two chambers:
\begin{align*}
B_+ & := \bR^{n-1} \times (|\epsilon|, +\infty),\\
B_- & := \bR^{n-1} \times (0,|\epsilon|)
\end{align*}
over which the Lagrangian torus fibers behave differently in a Floer-theoretic sense.

To explain what is going on, let us consider the genus 0 open Gromov-Witten invariants defined as in Definition \ref{defn:openGW}, or in other words, the virtual counts of {\em Maslov index 2} stable disks in $X$ bounded by fibers of $\rho$. Here, Maslov index 2 means, geometrically, that the stable disks intersect with the hypersurface $D_\epsilon$ at only one point with multiplicity one. As one moves from one chamber to the other by crossing the wall $H$, the virtual counts of Maslov index 2 disks bounded by the corresponding Lagrangian torus fiber would jump, exhibiting a {\em wall-crossing phenomenon}.

It was proved in \cite{CLL12} that fibers over the chamber $B_-$ bound only one (family) of Maslov index 2 disks, so that the corresponding generating function of open Gromov-Witten invariants has just one term with coefficient 1; while fibers over $B_+$ bound as many Maslov index 2 disks as a moment map Lagrangian torus fiber, so that the corresponding generating function of open Gromov-Witten invariants has (possibly infinitely) many terms. The resulting {\em wall-crossing factor} is exactly what we need in order to get the correct mirror.

More explicitly, the Lagrangian Floer superpotential $W^\text{LF}_q$ is given by
\begin{align*}
W^\text{LF}_q = \left\{
\begin{array}{ll}
z_0 \sum_{i=1}^{m} (1+\delta_i(q)) C_i z^{w_i} & \text{over $B_+$},\\
u & \text{over $B_-$},
\end{array}
\right.
\end{align*}
where
$$1 + \delta_i(q) = \sum_{\alpha\in H_2^{\text{eff}}(X;\bZ)} n_{\beta_i + d} q^d$$
is a generating function of genus 0 open Gromov-Witten invariants.
Here, both $z_0$ and $u$ denote the coordinate associated to the lattice point $(0,\ldots,0,1) \in \bZ^{n-1} \oplus \bZ$ (over $B_+$ and $B_-$ respectively), $z^w$ denotes the monomial $z_1^{w^1}\ldots z_{n-1}^{w^{n-1}}$ for $w=(w^1,\ldots,w^{n-1})\in\bZ^{n-1}$, $q^d$ denotes $\exp\left(-\int_d \omega_\bC\right)$ which can be expressed in terms of the complexified K\"ahler parameters $q_1, \ldots, q_r$, and $\beta_1, \ldots, \beta_m \in \pi_2(X,L)$ are the basic disk classes as before.

The semi-flat mirror $\check{X}_0 = TB_0/\Lambda$ in this case can be viewed as a gluing of the two complex charts
$$\check{X}_+ = TB_+/\Lambda, \quad \check{X}_- = TB_-/\Lambda$$
associated to the chambers $B_+$ and $B_-$ respectively. As aforementioned, nontrivial monodromy of the tropical affine structure on $B_0$ around the discriminant locus $\Gamma$ prevents the complex structure on $\check{X}$ from extending any further.

Here comes an important idea of Auroux \cite{Auroux07, Auroux09}: Floer theory says that $W^\text{LF}_q$ is an analytic function, so one should modify the gluing between the complex charts $\check{X}_+, \check{X}_-$ exactly by the wall-crossing factor:
$$\sum_{i=1}^{m} (1+\delta_i(q)) C_i z^{w_i}.$$
Such a modification would cancel the nontrivial monodromy of the complex structure around $\Gamma \subset B$ and gives an analytic superpotential on the corrected mirror $\check{X}$. We therefore arrive at the following theorem.\footnote{We remark that there is no scattering phenomenon in this case because there is only one wall, so no further corrections are needed.}

\begin{thm}[\cite{CLL12}, Theorem 4.37]
The SYZ mirror for the toric Calabi-Yau manifold $X$ is given by the family of affine hypersurfaces\footnote{Strictly speaking, this is the SYZ mirror for the complement $X\setminus D_\epsilon$ only; the mirror of $X$ itself is given by the Landau-Ginzburg model $(\check{X}, u)$.}
\begin{equation}\label{eqn:SYZmirror_toricCY}
\check{X}_q = \left\{(u,v,z_1,\ldots,z_{n-1}) \in \bC^2 \times (\bC^\times)^{n-1} : uv = \sum_{i=1}^{m} (1+\delta_i(q)) C_i z^{w_i} \right\}.
\end{equation}
\end{thm}

Notice that the SYZ mirror family \eqref{eqn:SYZmirror_toricCY} is entirely written in terms of complexified K\"ahler parameters and disk counting invariants of $X$, and it agrees, up to a mirror map, with the prediction \eqref{eqn:toricCY_mirror} using physical arguments \cite{Leung-Vafa98, CKYZ99, HIV00}.

\begin{example}\label{eg:SYZ_KP2}
The SYZ mirror of $X = K_{\bP^2}$ is given by
\begin{equation}\label{eqn:mirror_KP2}
\check{X} = \left\{(u,v,z_1,z_2)\in\bC^2\times(\bC^\times)^2 : uv = 1+\delta_0(q) + z_1 + z_2 + \frac{q}{z_1z_2} \right\},
\end{equation}
where $q$ is the K\"ahler parameter which measures the symplectic area of a projective line contained inside the zero section of $K_{\bP^2}$ over $\bP^2$, and
\begin{equation}\label{eqn:gen_fcn_openGW_KP2}
1+\delta_0(q) = \sum_{k=0}^\infty n_{\beta_0+kl} q^k
\end{equation}
is a generating function of genus 0 open Gromov-Witten invariants,
where $l\in H_2(K_{\bP^2};\bZ) = H_2(\bP^2;\bZ)$ is the hyperplane class in $\bP^2$, $q=\exp(-t)$ and $\beta_0$ is the basic disk class corresponding to the zero section $\bP^2 \subset K_{\bP^2}$, which is the only compact toric divisor in $K_{\bP^2}$.

By applying the open/closed formula in \cite{Chan11} and a flop formula for closed Gromov-Witten invariants \cite{Li-Ruan01}, the genus 0 open Gromov-Witten invariants $n_{\beta_0+kl}$ can be expressed in terms of {\em local BPS invariants} of $K_{\mathbb{F}_1}$, where $\mathbb{F}_1$ is the blowup of $\bP^2$ at one point \cite{LLW11}. More specifically, if
$$e, f \in H_2(K_{\mathbb{F}_1};\bZ) = H_2(\mathbb{F}_1;\bZ)$$
are the classes represented by the exceptional divisor and fiber of the blowup $\mathbb{F}_1 \to \bP^2$ respectively, then
$$n_{\beta_0+kl} = N_{0,0,kf+(k-1)e}^{K_{\mathbb{F}_1}},$$
where the right-hand side is the so-called {\em local BPS invariant} for the class $kf+(k-1)e\in H_2(K_{\mathbb{F}_1},\bZ)$. The latter invariants have been computed by Chiang, Klemm, Yau and Zaslow and listed in the ``sup-diagonal'' of Table 10 on p. 56 in \cite{CKYZ99}:
\begin{align*}
n_{\beta_0+l}  & = -2,\\
n_{\beta_0+2l} & = 5,\\
n_{\beta_0+3l} & = -32,\\
n_{\beta_0+4l} & = 286,\\
n_{\beta_0+5l} & = -3038,\\
n_{\beta_0+6l} & = 35870,\\
               &\vdots
\end{align*}
Hence, the first few terms of the generating function $1+\delta_0(q)$ is given by
$$1+\delta_0(q) = 1 - 2q + 5q^2 - 32q^3 + 286q^4 - 3038q^5 + \cdots,$$
which, when compared with the power series \eqref{eqn:KP2_integral_series}, suggests that the coefficients of $f(q)$ are exactly the virtual counts $n_{\beta_0 + kl}$ of holomorphic disks. This is indeed the case as we will see in the next section.
\end{example}

We remark that the SYZ construction can be carried out also in the reverse direction; see the work of Abouzaid-Auroux-Katzarkov \cite{AAK12}.
For example, starting with the conic bundle \eqref{eqn:mirror_KP2}, it is possible to construct a Lagrangian torus fibration using techniques similar to that in \cite{Goldstein01, Gross01}. Although the discriminant locus in this case is of real codimension {\em one} (which is an {\em amoeba} and a thickening of the graph $\Gamma$ in Figure \ref{KP2_base}), one can still construct the SYZ mirror by essentially the same procedure and get back the toric Calabi-Yau 3-fold $K_{\bP^2}$.\footnote{Strictly speaking, the SYZ mirror of \eqref{eqn:mirror_KP2} is the complement of a smooth hypersurface (namely $D_\epsilon$) in $K_{\bP^2}$.} This demonstrates the involutive property of the SYZ mirror construction.


\subsection{Period integrals and disk counting -- a conjecture of Gross-Siebert}

As have been anticipated by Gross and Siebert in their spectacular program \cite{Gross-Siebert03, Gross-Siebert-logI, Gross-Siebert-logII, Gross-Siebert-reconstruction}, the SYZ construction should give us a mirror family \eqref{eqn:SYZmirror_toricCY} which is inherently written in {\em canonical coordinates}; cf. also \cite{Kontsevich-Soibelman06}. This is an indication why the SYZ construction is the correct and natural construction. Just like the open mirror theorem for compact toric semi-Fano manifolds (Section \ref{sec:open_mirror_thm}), this expectation can be reformulated in a more succinct way via the {\em SYZ map}

\begin{defn}\label{defn:SYZ_map_toricCY}
We define the {\em SYZ map}
\begin{align*}
\phi:\mathcal{M}_A & \to \check{\mathcal{M}}_B,\\
q=(q_1,\ldots,q_r) & \mapsto \phi(q)=(t_1(q),\ldots,t_r(q)),
\end{align*}
by
\begin{align*}
t_a(q) = q_a \prod_{i=1}^{m} \left(1 + \delta_i(q)\right)^{D_i \cdot \gamma_a},\quad a=1,\ldots,r.
\end{align*}
\end{defn}

Then we have the following conjecture.
\begin{conj}[\cite{CLL12}, Conjecture 1.1]\label{conj:SYZ_map=period}
There exist integral cycles $\Gamma_1,\ldots,\Gamma_r$, forming part of an integral basis of the middle homology group $H_n(\check{X}_t;\bZ)$, such that
\begin{align*}
q_a = \exp\left(-\int_{\Gamma_a}\check{\Omega}_{\phi(q)}\right)\text{ for $a=1,\ldots,r$,}
\end{align*}
where $\phi(q)$ is the SYZ map defined in Definition \ref{defn:SYZ_map_toricCY} in terms of generating functions $1+\delta_i(q)$ of the genus 0 open Gromov-Witten invariants $n_{\beta_i+\alpha}$. In other words, the SYZ map is inverse to a mirror map.
\end{conj}

Conjecture \ref{conj:SYZ_map=period} not only provides an enumerative meaning to the inverse of a mirror map which is defined by period integrals, but also explains the mysterious integrality property of the coefficients of its Taylor series expansions which was observed earlier, e.g. by Zhou \cite{Zhou12}.

This conjecture can be regarded as a more precise reformulation of a conjecture of Gross and Siebert \cite[Conjecture 0.2 and Remark 5.1]{Gross-Siebert-reconstruction} where they predicted that period integrals of the mirror are intimately related to (virtual) counting of tropical disks (instead of holomorphic disks) emanating from the singularities in the base of a (special) Lagrangian torus fibration on a {\em compact} Calabi-Yau manifold; see \cite[Example 5.2]{Gross-Siebert11a} where Gross and Siebert observed a relation between the so-called {\em slab functions}, which play a key role in their program, and period integral calculations for $K_{\bP^2}$ in \cite{Graber-Zaslow02}.

Numerical evidences of Conjecture \ref{conj:SYZ_map=period} were first given for the toric Calabi-Yau surface $K_{\bP^1}$ and the toric Calabi-Yau 3-folds $\mathcal{O}_{\bP^1}(-1)\oplus\mathcal{O}_{\bP^1}(-1)$, $K_{\bP^2}$ (as we have seen in Example \ref{eg:SYZ_KP2}), $K_{\bP^1\times \bP^1}$ in \cite[Section 5.3]{CLL12}. In \cite{LLW12}, Conjecture \ref{conj:SYZ_map=period} was proved for all toric Calabi-Yau surfaces (which are all crepant resolutions of $A_n$-surface singularities).

In \cite{CLT11}, a slightly weaker version of Conjecture \ref{conj:SYZ_map=period}, where the cycles $\Gamma_1,\ldots,\Gamma_r$ are allowed to have complex coefficients, was proved for the case when $X$ is the total space of the canonical line bundle $K_Y$ over a toric Fano manifold $Y$. The proof is based on the open/closed equality \eqref{eqn:open_closed} proved in \cite{Chan11} and the mirror theorem of Givental \cite{Givental98} and Lian-Liu-Yau \cite{LLY-III}. In particular, for $X = K_{\bP^2}$, this verifies the claim that the generating function $1+\delta_0(q)$ of genus 0 open Gromov-Witten invariants in \eqref{eqn:gen_fcn_openGW_KP2} is equal to the function $f(q)$ in \eqref{eqn:KP2_integral_series} computed using period integrals.

Finally, the following theorem was proved in \cite{CCLT13}:
\begin{thm}[\cite{CCLT13}, Corollary 1.7]\label{thm:openGW_period}
Given any toric Calabi-Yau manifold $X$, there exists a collection $\{\Gamma_1,\ldots,\Gamma_r\} \subset H_n(\check{X};\bC)$ of linearly independent cycles in the SYZ mirror $\check{X}$ (defined in \eqref{eqn:SYZmirror_toricCY}) such that
\begin{equation*}
q_a = \exp\left(-\int_{\Gamma_a}\check{\Omega}_{\phi(q)}\right), \quad a = 1,\ldots,r,
\end{equation*}
where $q_a$'s are the complexified K\"ahler parameters in the complexified K\"ahler moduli space $\mathcal{M}_A$ of $X$, and $\phi(q)$ is the SYZ map in Definition \ref{defn:SYZ_map_toricCY}, defined in terms of the generating functions $1+\delta_i(q)$ of genus 0 open Gromov-Witten invariants $n_{\beta_i+\alpha}$.
\end{thm}

To prove this theorem, we employ orbifold techniques. Indeed, the SYZ mirror construction and the above relation between period integrals and open Gromov-Witten theory have natural extensions to toric Calabi-Yau {\em orbifolds}. This is an important point since it enables us to formulate and prove an {\em open} analogue of the {\em crepant transformation conjecture} for toric Calabi-Yau varieties; we refer the reader to \cite[Section 1]{CCLT13} for an introduction to these results.

\begin{remark}
In order to prove the original Conjecture \ref{conj:SYZ_map=period}, what we lack is a construction of {\em integral} cycles whose periods have specific logarithmic terms.

Indeed, such cycles have been constructed by Doran and Kerr in \cite[Section 5.3 and Theorem 5.1]{Doran-Kerr11} when $X$ is the total space of the canonical line bundle $K_S$ over a toric del Pezzo surface $S$. It should not be difficult to extend their construction to higher dimensions and even more general toric Calabi-Yau varieties; cf. the discussion in \cite[Section 4]{Doran-Kerr13}. We thank Chunk Doran for useful discussions related to this problem.

On the other hand, in a recent paper \cite{Ruddat-Siebert14}, Ruddat and Siebert gave another construction of such integral cycles by applying tropical methods. Though they worked only in the compact Calabi-Yau case, it is plausible that their method can be generalized to handle the toric Calabi-Yau case as well. We thank Helge Ruddat for pointing this out, who mentioned that they were planning to add the toric Calabi-Yau case to their paper \cite{Ruddat-Siebert14} in the future.
\end{remark}

\begin{remark}
In their ICM lecture \cite{Gross-Siebert_ICM}, Gross and Siebert sketched an alternative proof of Theorem \ref{thm:openGW_period}, with counting of holomorphic stable disks replaced by counting of tropical disks (or broken lines), via {\em logarithmic Gromov-Witten theory}. It is therefore natural to expect that our SYZ mirror family \eqref{eqn:SYZmirror_toricCY} coincides with the mirror constructed from the Gross-Siebert program. Very recently, this was confirmed by Lau \cite{Lau14} who showed that the generating functions $1 + \delta_i(q)$ are indeed equal to the slab functions of Gross and Siebert.
\end{remark}

\section*{Acknowledgment}
The author would like to thank Cheol-Hyun Cho, Siu-Cheong Lau, Naichung Conan Leung and Hsian-Hua Tseng for very fruitful collaborations which lead to much of the work described here and the referees for valuable comments and pointing out a number of errors.

Some of the research described in this paper were substantially supported by a grant from the Research Grants Council of the Hong Kong Special Administrative Region, China (Project No. CUHK404412).

\bibliographystyle{amsplain}
\bibliography{geometry}

\end{document}